\documentclass[final,1p,times]{elsarticle}

\usepackage{amssymb,amsmath,bbm,amsthm,dsfont,amsthm, color}

\usepackage{lineno}

\journal{arXiv}

\begin{document}

\def\Xint#1{\mathchoice
   {\XXint\displaystyle\textstyle{#1}}%
   {\XXint\textstyle\scriptstyle{#1}}%
   {\XXint\scriptstyle\scriptscriptstyle{#1}}%
   {\XXint\scriptscriptstyle\scriptscriptstyle{#1}}%
   \!\int}
\def\XXint#1#2#3{{\setbox0=\hbox{$#1{#2#3}{\int}$}
     \vcenter{\hbox{$#2#3$}}\kern-.5\wd0}}
\def\ddashint{\Xint=}
\def\dashint{\Xint-}
\newcommand*{\id}{\mathds{1}}
\newcommand*{\zahlen}{\mathbb{Z}}
\newcommand*{\en}{\mathbb{N}}
\newcommand*{\er}{\mathbb{R}}
\newcommand{\esssup}{\mathop{\text{ess\:sup}}}
\newcommand{\htimes}{\mathop{\text{\large$◊$}}}
\newcommand{\hexists}{\mathop{\text{\LARGE$\exists$}}}
\newcommand{\hforall}{\mathop{\text{\LARGE$\forall$}}}

\newtheorem{l1}{Lemma}
\newtheorem{Tw}{Theorem}
\newtheorem*{Tw*}{Theorem}
\newtheorem{cor}{Corollary}
\newtheorem{mydef}{Definition}

\begin{frontmatter}



\title{$L^\infty$ a priori bounds for gradients of solutions to quasilinear inhomogenous fast-growing parabolic systems}


\author{Jan Burczak\fnref{label2}}
 \fntext[label2]{The author is a Ph.D. student in the International Ph.D. Projects Programme of Foundation for Polish Science within the Innovative Economy Operational Programme 2007-2013 (Ph.D. Programme: Mathematical Methods in Natural Sciences). This publication is partly supported by the National Science Centre (NCN) grant
no. 2011/01/N/ST1/05411.}
\address{Institute of Mathematics, Polish Academy of Sciences, ul. \'Sniadeckich 8,
00-956~Warszawa}
\ead{jb@impan.pl}

\begin{abstract}
We prove boundedness of gradients of solutions to quasilinear parabolic system, the main part of which is a generalization to p-Laplacian and its right hand side's growth depending on gradient is not slower (and generally strictly faster) than $p-1$. This result may be seen as a generalization to the classical notion of a controllable growth of right hand side, introduced by Campanato, over gradients of p-Laplacian-like systems. Energy estimates and nonlinear iteration procedure of a~Moser type are cornerstones of the used method. 

\end{abstract}


\end{frontmatter}


\section{Introduction}
\noindent

\subsection{General statement of the problem.}
\noindent
We are interested in obtaining a local boundedness of gradients of solutions to a following parabolic system in $\Omega \subset \er^n$
\begin{equation*}
u^i_t - (A^i_{\alpha} ( \nabla u))_{x_\alpha} = f^i (x,t, \nabla u)  \qquad i = 1,..,N
\end{equation*}
where the main part is a generalization of p-Laplacian and the right-hand-side grows as $1+ |\nabla u|^w$ or $ |\nabla u|^w$ with $w$ specified further. We say that a right-hand-side is a fast-growing one, when $w > p-1$ holds. \newline
The existing literature on the regularity issue of parabolic equations and systems is impressive. Let us recall that for equations the existing results are quite strong: even for the  right-hand-side growth of $1+|\nabla u|^p$ one obtains $C^{1, \alpha}$ regularity of solutions: see classic monograph \citet{S-L-U68} for the case $p=2$ and \citet{diBbook93} for $p \in (1, \infty)$. Many further generalizations are possible: for instance in \citet{Bar-Sou04} the right-hand-side takes the form $e^{u} |\nabla u|^p$, which suffices for a boundedness of $\nabla u$. Moreover, this growth condition seems to be optimal, because there are blowup results for gradients of solutions to equations, the right-hand-sides of which grow faster than $p$ - compare \citet{Sou02}. In case of systems, the regularity results are much weaker. One can construct irregular (i.e. unbounded or discontinous) functions, which solve homogenous parabolic systems. For $n>2$ it suffices for irregularity that the coefficients $A(x,t)$ of the main part are discontinous (and still bounded) or that there is a relevant non-diagonality of the main part - for details, consult \citet{Arkh08}. Nevertheless, there are many classes of main parts, which allow for higher regularity (even $C^{1, \alpha}$) in the homogenous case; these are: having structure close to Laplacian or p-Laplacian, like those studied in \citet{S-L-U68} or \citet{diBbook93}, respectively, or having main part depending solely on $\nabla u$: see well-known paper by \citet{Nec-Sve91} or more extensive research done by \citet{Bae-Choe04}. As these papers consider homogenous systems, one may ask a natural question: what inhomogenous counterparts of such systems remain, in a certain sense, regular? The general answer is unknown, but there are several hints: on one hand, for the right-hand-side growing like $1+ |\nabla u|^{p-1}$ the regularity of the homogenous case seems to be retained - see \citet{diBbook93}; on the other, unlike for equations, one cannot have the right-hand-side growing as fast as $ |\nabla u|^{p}$ without further assumptions, even in case of a system with the simplest main part, i.e. an inhomogenous heat system. Recall the classical counterexample: for $n \ge 3$ bounded but discontinous function $u(x) = \frac{x}{|x|}$ with unbounded weak derivatives  solves $u^i_t  - \Delta u^i = u^i |\nabla u|^2 \quad \left[ = (n-1) \frac{x^i}{|x|^3}\right]$, for details - see \citet{Str84}. It turns out that in case of an inhomogenous system for $p=2$ one has to additionally assume a certain smallness in order to obtain regularity - for details, refer to \citet{Tolks86, Pin07, Idone87} or even the classical \citet{S-L-U68}.
The regularity issue for a general nonlinear inhomogenous parabolic system with the right-hand-side growing at the rate $1+|\nabla u|^w$ for $w$ possibly close to $p$,  homogenous counterparts of which enjoy regularity, is not fully researched, especially for the case $p \neq 2$. There are several approaches to answering this question: some authors relax the notion of regularity by resorting to partial regularity - see for example classical papers of Italian school: \citet{Cam84, Gia-Str82} and newer ones: \citet{ Fan97, Spe-Freh09,  Mis02, Duz-Ming05}; or by demanding a high integrability-type regularity\footnote{such results are especially interesting, as our result may be easily strenghtened via higher regularity}, like in \citet{Nau90, Kin-Lew00} or \citet{Bens-Freh09} (in the last paper the growth of the right-hand-side may be polynomially arbitrarily large!). Certain systems with peculiar structure or two-dimensional ones (or at least close to them in some sense) enjoy also high regularity, even if they are much more general than a Stokes-type system; for results in this direction compare papers of Seregin, Arkhipova, Frehse, Kaplicky (and many others), for instance: \citet{Arkh08, Nau-Wolf-Wolff98, Kap05, Zaj-Ser07}. \newline
In this note we focus on deriving a full regularity result, more precisely: the local boundedness of gradients, for a class of quasilinliear parabolic inhomogenous systems. Our goal is twofold: firstly to obtain results for a general inhomogenous parabolic system, the main part of which is analogous to the system considered in \citet{Bae-Choe04}, while retaining possibly general growth conditions for the right-hand-side. Secondly, to sharpen these results with respect to growth of the right-hand-side, restricting ourselves to less general systems, being close to $p-$Laplacian. For similar result on the level of solutions, compare \citet{Gio-O'Le04}. \newline 
Let us emphasize that we proceed in a manner typical for the regularity approach: we assume existence of solution $u$ in a given class, which is often a deep problem itself, from which we derive higher regularity. Moreover, we concentrate on a priori estimates while conducting the proofs: the rigorous version of computations is commented on in the conclusion. 

\subsection{General definitions and assumptions.}
\noindent
Consider parabolic problem in $\Omega \subset \er^n$
\begin{equation}
u^i_t - (A^i_{\alpha} ( \nabla u))_{x_\alpha} = f^i (x,t, \nabla u)  \qquad i = 1,..,N  \label{i1}
\end{equation}
As all our results have a local character, any further specification of $\Omega$ is irrelevant. \newline
We say that a vector valued function $u \in L^\infty (0,T;L^2(\Omega)) \cap L^p(0,T;W^{1,p}(\Omega))$ is a weak solution to \eqref{i1} iff
\begin{equation*}
\int_{\Omega_T} - u^i \phi^i_t + A^i_{\alpha} ( \nabla u)  \phi^i_{x_\alpha} dx dt= \int_{\Omega_T} f^i (x,t, \nabla u) \phi^i dxdt\qquad \hforall_{\phi \in C^\infty_0 (\Omega_T)} 
\end{equation*}
Globally, the following notions will be used:
\begin{itemize}
\item $\delta^\alpha_\beta$ denotes the Kronecker delta,
\item $Q_R (x_0 , t_0 )$ denotes parabolic cylinder, i.e. $ B_R (x_0 ) \times ( t_0 - R^p, t_0)$; when possible, cut short to $Q_R$,
\item $\eta_{\rho, R} \in C_0^\infty (Q_R)$ denotes standard parabolic cutoff function for $ Q_\rho \subset Q_R$, when possible, cut short to $\eta$, which satisfies: $\eta = 1$ in $Q_\rho$, $\eta = 0$  outside $Q_R$, $|\nabla \eta| \le c (R - \rho)^{-1}$, $ |\eta_{,t}| \le c (R - \rho)^{-p}$.
\end{itemize}
Throughout the article, summation over repeated indices is in use.

\subsection{The structure of results.}
We show our results in the following order:
\begin{enumerate}
\item First, we derive a general result for inhomogenous version of system analyzed in \citet{Bae-Choe04}, where ellipticity assumptions for the main part are generalized by introducing exponent $q$  (Theorem \ref{thm2}). Here, loosely speaking, admissible growth for the right-hand-side is $1 + |\nabla u|^{p-1}$, so this result may be seen as parallel to \citet{diBbook93}.
\item Next, we admit faster growths for the right-hand-side, at the cost of assuming that the main part is closer to p-Laplacian, in the sense that it is not enriched with terms involving $q>p$ (Theorem \ref{thm2_p=q}). 
\item Finally, we state the result for the least general case, i.e. for the 3D p-Laplacian with the right-hand-side growing as $| \nabla u|^w$ (Theorem \ref{thm_intr}). 
\end{enumerate}
Since the last result seems to be the most traceable one,  let us give incentive to studying the technical remainder of this paper by stating Theorem \ref{thm_intr} now: \newline
\noindent
Consider $u = (u^i) \in L^\infty (0,T; L^2(\Omega)) \cap L^p(0,T; W^{1,p}(\Omega)) $ solving the $3$D p-Laplace system
\begin{equation*}
u^i_t - {\rm div} (| \nabla u|^{p-2} \nabla u^i ) = f^i (t, x, u, \nabla u) \quad i=1,2,3
\label{wst}
\end{equation*}

\begin{Tw}
Let $\Omega$ be an arbitrary domain. Assume a growth condition: $|f^i  (x,t, \nabla u)| \le | \nabla u |^w, \; w\le p$ and initial integrability\footnote{From existence one has $  {\tilde p} = p$, so this assumption may be void. It only helps to quantify the results when we have some additional knowledge on integrability.} $ |\nabla u |_{L^{\tilde p}_{loc}} < \infty$. 
If \textbf{one} of the following conditions is fullfilled:
\begin{enumerate}
\item $w \le p-1, \: {\tilde p} = p$;
\item $w \in \left[\frac{p}{2}, \frac{ {\tilde p} + 4p - 3}{5} \right)$;
\item $w \le  \frac{ p }{2}$ and $p > 2 - \frac{2}{3}  {\tilde p} $
\end{enumerate}
then $\nabla u$ is bounded:
\begin{equation*}
|\nabla u |_{L^{\infty}_{loc}} < C|\nabla u |_{L^{\tilde p}_{loc}} 
\end{equation*}
\label{thm_intr}
\end{Tw}
\noindent
For the proof, see the end of the next section. \newline
Observe that for $w > p-1$ in the degenerate case (i.e. $p \ge 2$) we have $w \ge  \frac{ p }{2}$, so point $2.$ applies. In such case merely from existence, i.e for ${\tilde p} = p$,  one has $w < p - \frac{3}{5}$. This corresponds for $p=2$ with Campanato's notion of controllable growth, stating that to obtain regularity, the growth of right-hand-side must be smaller than $\frac{7}{5}$  (see for example \citet{Cam84}). Therefore, this result may be seen as a generalization of  classical results over gradients of p-Laplacian-like systems. Utilizing results on high integrability of certain systems, one may relax the growth condition further. For example, for the system analysed in \citet{Nau90} we have $w < \varepsilon + \frac{9}{5}$, because $\nabla u \in L^{ \varepsilon + 4} (\Omega_T)$, which can be taken as the initial integrability $L^{\tilde p} $ .
\newline
In all our  theorems there is no explicit assumption that $w < p$. In fact, the inequality $w \le p$ is enforced by the rigorous treatment of energy estimates. Simultaneously we know from the counterexample recalled in the introduction, that $w = p$  is generally not admissible. Therefore, our results can be viewed as a way to quantify the possible boundedness of gradients by means of higher integrability. As in the just mentioned case of Theorem \ref{thm_intr}, for $p=2$ one has $w < p - \frac{3}{5}$, which can be boosted in some cases to $w < \varepsilon + \frac{9}{5}$, because $\nabla u \in L^{ \varepsilon + 4}$. For $w=p$ one would need $\nabla u \in L^{ \varepsilon + 5}$.

\section{Boundedness of gradient of solution}
\noindent
As outlined in the introduction, first we prove the general theorem. As the cornerstone of the analysis is the energy method, we derive formal estimates for the sake of transparency. For  a rigorous justification of the formal estimates please consult the conclusion of this note.
We analyze solutions of
\begin{equation}
u^i_t - (A^i_{\alpha} ( \nabla u))_{x_\alpha} = f^i (x,t, \nabla u)  \qquad i = 1,..,N  \label{bg1}
\end{equation}
where the main part comes from \citet{Bae-Choe04} and the right-hand-side grows as $1 + |\nabla u|^w, \; w \le \delta $ for a certain $ \delta \le p$, obtaining the boundedness of $\nabla u$. More precisely, one has the following:
\begin{Tw}
Under the following assumptions:
\begin{enumerate}
\item[(A0)] ellipticity-type: $A^i_{\alpha}$ is given by potential $F \in C^2 (\er), \: F' (0) \ge 0,$ as follows
\begin{equation*}
A^i_\alpha (Q) = (F(|Q|)_{Q^i_\alpha}
\end{equation*}
and $F$ enjoys ellipticity
\begin{equation*}
(F(|Q|)_{Q^i_\alpha Q^j_\beta} \zeta^i_\alpha \zeta^j_\beta  \in [ \lambda  | Q|^{p-2} ,  \lambda^{-1} ( | Q|^{p-2}+   |Q|^{q-2} )] |\zeta|^2
\end{equation*}
where:
\begin{equation*}
1 < p \le q < p+1 < \infty
\end{equation*}
\item[(A1)] growth-type
\begin{equation*}
|f^i| \le + c_1 | \nabla u |^w + c_2
\label{bg_growth}
\end{equation*}
where
\begin{equation*}
p \ge w \ge 0, \quad c_i \in L^\infty (\Omega_T)
\end{equation*}
\item[(A2)] initial integrability
\begin{equation*}
\nabla u \in L_{loc}^{s_0 + M} (\Omega_T), \quad M := \max (2, p, 2q-p, w+1, 2w-p+2)  
\end{equation*}
with $s_0 $ satisfying
\begin{equation*}
s_0 \ge 0, \qquad
s_0 + 2 +\frac{np}{2} - \frac{M n}{2}  > 0 
\end{equation*}
and
\begin{equation*}
\begin{cases}
s_0  > p-2 \; for\;  c_2 \neq 0 \\
s_0  > p-2w-2 \; for\;  c_2 = 0
\end{cases}
\label{ass_alt1}
\end{equation*}
\end{enumerate}
gradient of solution to \eqref{bg1} is locally bounded;
moreover, for any $Q_{R_{0}} \subset \Omega$ with $R_{0} < 1$ the following inequality holds
\begin{equation*}
 |\nabla u|_{L^{\infty} (Q_{\frac{R_0}{2}})} \le C \left( \dashint_{Q_{R_{0}}} |\nabla u|^{s_0 +M} dxdt \right)^{\frac{1}{s_0 +2 + \frac{np}{2} - \frac{Mn}{2}}} + C
 \end{equation*}
\label{thm2}
\end{Tw}

\begin{proof}
First we derive formal energy inequalities, then we implement an iteration scheme. \newline
\noindent
Differentiate formally system \eqref{bg1} to obtain
\begin{equation}
u^i_{t x_\gamma} - (A^i_{\alpha, \:  u^j_{x_\beta}} (\nabla u)   u^j_{x_\beta  x_\gamma}  )_{x_\alpha}=( f^i (x,t, \nabla u) )_{x_\gamma} \qquad i = 1,..,N  \label{bg_p1}
\end{equation}
testing \eqref{bg_p1} by $u^i_{x_\gamma} | \nabla u|^s \eta^2$ one gets

\begin{multline}
\left[ \frac{1}{s+2} \frac{d}{dt} \int_{B_R} |\nabla u|^{s+2} \eta^2   dx \right]+ 
\underbrace{\left[  \int_{Q_R} A^i_{\alpha \: , \:  u^j_{x_\beta}} (\nabla u)   u^j_{x_\beta  x_\gamma} ( |\nabla u|^s u^i_{x_\gamma x_\alpha} \eta^2 + s |\nabla u|^{s-2} u^i_{x_\gamma} u^k_{x_\alpha x_\delta}   u^k_{x_\delta} \eta^2 ) dx \right]}_{I} = \\
 - \int_{B_R} f^i (x,t, \nabla u) [ u^i_{x_\gamma x_\gamma}  | \nabla u|^s \eta^2 + s u^i_{x_\gamma}  | \nabla u|^{s-2} u^k_{x_\delta x_\gamma}  u^k_{x_\delta} \eta^2  + 
2 u^i_{x_\gamma}  | \nabla u|^s \eta \eta_{x_\gamma}  ]  dx \\
 + \frac{2}{s+2} \int_{B_R} |\nabla u|^{s+2} \eta \eta_{t} dx -  \int_{B_R}  A^i_{\alpha , \:  u^j_{x_\beta}} (\nabla u) u^k_{x_\beta x_\gamma}  u^k_{x_\gamma} |\nabla u|^{s}  \eta \eta_{x_\alpha} dx \label{bg_p2}
\end{multline}
Consider $I$. Utilizing ellipticity assumption $(A0)$ with $\zeta^l_\rho := u^l_{x_\beta  x_\rho }$ one estimates the first summand of $I$ as follows
\begin{equation}
\int_{B_R} A^i_{\alpha \: '(\nabla u)^j_\beta} (\nabla u)   u^j_{x_\beta  x_\gamma}  |\nabla u|^s u^i_{x_\gamma x_\alpha} \eta^2 dx \ge \lambda \int_{B_R} |\nabla u|^{p-2}   |\nabla^2 u|^2  |\nabla u|^{s}  \eta^2 dx
 \label{bg_p3}
\end{equation}
 and because $A^i_\alpha$ is given by potential $F$, from differentiation we estimate the second summand of $I$
 \begin{multline}
 \int_{B_R} A^i_{\alpha \: '(\nabla u)^j_\beta} (\nabla u)   u^j_{x_\beta  x_\gamma} s |\nabla u|^{s-2} u^i_{x_\gamma} u^k_{x_\alpha x_\delta}   u^k_{x_\delta} \eta^2 dx =\\  
 \int_{B_R} \left[ F''(|\nabla u |) \frac{u^i_{x_\alpha} u^j_{x_\beta}}{|\nabla u |^2}+ F'(|\nabla u |) \left( \frac{ \delta^i_{x_\alpha} \delta^j_{x_\beta}}{|\nabla u |} -  \frac{u^i_{x_\alpha} u^j_{x_\beta}}{|\nabla u |^3} \right)   \right] u^j_{x_\beta  x_\gamma} s |\nabla u|^{s-2} u^i_{x_\gamma} u^k_{x_\alpha x_\delta}   u^k_{x_\delta} \eta^2 dx \\
  =  s \int_{B_R}  \eta^2 F''(|\nabla u |) {|\nabla u |^{s-4}} \underbrace{( u^i_{x_\gamma}  u^j_{x_\beta} u^j_{x_\beta x_\gamma} )}_{c^i} \underbrace{ (u^i_{x_\alpha}  u^k_{x_\delta} u^k_{x_\delta x_\alpha})}_{c^i} \ dx  +\\ 
 s \int_{B_R}  \eta^2 F'(|\nabla u |)  {|\nabla u |^{s-5}} ( \underbrace{u^k_{x_\delta} u^k_{x_\delta x_\alpha}  u^i_{x_\gamma} u^i_{x_\gamma x_\alpha} |\nabla u |^2}_{= \frac{1}{4} (|\nabla u|^2)_{x_\alpha}  (|\nabla u|^2)_{x_\alpha}  | \nabla u|^2 } - \underbrace{u^i_{x_\alpha}  u^i_{x_\gamma}  u^j_{x_\beta}  u^j_{x_\beta x_\gamma}  u^k_{x_\delta}  u^k_{x_\delta x_\alpha}}_{ = \frac{1}{4} u^i_{x_\alpha}  (|\nabla u|^2)_{x_\alpha} u^i_{x_\gamma}  (|\nabla u|^2)_{x_\gamma}  } ) dx
  \label{bg_p4}
\end{multline}
From the ellipticity assumption $(A0)$ one has $F''(|s|) \ge 0, \: F'(0) \ge 0$ therefore it holds: $ F' (|s|) \ge 0$. This, in conjunction with the following computation: $ u^i_{x_\alpha}  (|\nabla u|^2)_{x_\alpha} u^i_{x_\gamma}  (|\nabla u|^2)_{x_\gamma} \le  | \nabla u |^2 | \nabla |\nabla u|^2 |^2 =   | \nabla u |^2  (|\nabla u|^2)_{x_\beta}  (|\nabla u|^2)_{x_\beta}  $, implies that equation \eqref{bg_p4} takes the form
 \begin{equation}
 \int_{B_R} A^i_{\alpha \: '(\nabla u)^j_\beta} (\nabla u)   u^j_{x_\beta  x_\gamma} s |\nabla u|^{s-2} u^i_{x_\gamma} u^k_{x_\alpha x_\delta}   u^k_{x_\delta} \eta^2 dx \ge 0
  \label{bg_p5}
\end{equation}

Summing up \eqref{bg_p3} and \eqref{bg_p5} we conclude that $I$ satisfies
\begin{multline}
 \int_{B_R} A^i_{\alpha \: '(\nabla u)^j_\beta} (\nabla u)   u^j_{x_\beta  x_\gamma} ( |\nabla u|^s u^i_{x_\gamma x_\alpha} \eta^2 + s |\nabla u|^{s-2} u^i_{x_\gamma} u^k_{x_\gamma x_\delta}   u^k_{x_\delta} \eta^2 ) dx \\
 \ge  \lambda \int_{B_R} |\nabla u|^{p-2}   |\nabla^2 u|^2  |\nabla u|^{s}  \eta^2 dx
   \label{bg_p6}
\end{multline}
Inputting inequality  \eqref{bg_p6} into   \eqref{bg_p2} we arrive at
\begin{multline}
\frac{1}{s+2} \frac{d}{dt} \int_{B_R} |\nabla u|^{s+2} \eta^2 dx +  \lambda \int_{B_R} |\nabla u|^{p-2}   |\nabla^2 u|^2  |\nabla u|^{s}  \eta^2 dx \le \\ \int_{B_R} f^i (x,t, \nabla u) [ u^i_{x_\gamma x_\gamma}  | \nabla u|^s \eta^2 + s u^i_{x_\gamma}  | \nabla u|^{s-2} u^k_{x_\delta x_\gamma}  u^k_{x_\delta} \eta^2  + 
2 u^i_{x_\gamma}  | \nabla u|^s \eta \eta_{x_\gamma}  ] dx \\
 + \frac{2}{s+2} \int_{B_R} |\nabla u|^{s+2} \eta \eta_{t} dx -  \int_{B_R}  A^i_{\alpha (\nabla u)^j_\beta} (\nabla u) u^k_{x_\beta x_\gamma}  u^k_{x_\gamma} |\nabla u|^{s}  \eta \eta_{x_\alpha} dx \\ \le
c \int_{B_R} [ 1 + | \nabla u|^w ] [ (1+ s) |\nabla^2 u|  | \nabla u|^s \eta^2 + 
2 | \nabla u|^{s+1} \eta | \nabla \eta| ] dx + \\
 \frac{c}{s+2} \int_{B_R} |\nabla u|^{s+2} \eta |\eta_{t}| dx + \frac{c}{\lambda}  \int_{B_R} [| \nabla u|^{p-2} + | \nabla u|^{q-2}]  |\nabla^2 u|  | \nabla u|^{s+1}   \eta | \nabla \eta| dx
  \label{bg_p7}
\end{multline}
where the last inequality is valid in view of growth $(A1)$ and ellipticity   $(A0)$  assumptions. \\
Absorb $ |\nabla^2 u|$ from the right-hand-side of   \eqref{bg_p7} using Young's inequality and integrate with respect to time
\begin{multline}
\frac{1}{s+2} \sup_t \int_{B_R} |\nabla u|^{s+2} \eta^2 dx +  (\lambda - \varepsilon) \int_{Q_R} |\nabla u|^{p-2}   |\nabla^2 u|^2  |\nabla u|^{s}  \eta^2 dxdt \le \\
c \int_{Q_R} ( 1 + | \nabla u|^w )  | \nabla u|^{s+1} \eta | \nabla \eta | dxdt
 + \frac{c}{s+2} \int_{Q_R} |\nabla u|^{s+2} \eta |\eta_{t}|  dxdt+ \\
 c \int_{Q_R} | \nabla u|^s  [| \nabla u|^{p} + | \nabla u|^{2q-p}]   \eta^2 | \nabla \eta|^2 dxdt +
 c(1+s)  \int_{Q_R} \eta^2 | \nabla u|^s [| \nabla u|^{2w - p +2} + | \nabla u|^{2-p}]   dxdt
  \label{bg_p9}
\end{multline}
By estimates for derivatives of the cutoff function $\eta$ we obtain
\begin{multline}
\frac{1}{s+2} \sup_t \int_{B_R} |\nabla u|^{s+2} \eta^2 dx + (\lambda - \varepsilon) \int_{Q_R} |\nabla u|^{p-2}   |\nabla^2 u|^2  |\nabla u|^{s}  \eta^2  dxdt \le \\
 \frac{c}{(R -\rho)^{\max(2,p)} }\int_{Q_R} | \nabla u|^s  [ | \nabla u|^{p} + | \nabla u|+ | \nabla u|^{w+1} + \\
 (2+s)^{-1}  | \nabla u|^2 + | \nabla u|^{p} + | \nabla u|^{2q-p} + (1+s) ( | \nabla u|^{2w - p +2} + | \nabla u|^{2-p} ) ] dxdt
 \label{bg_p10}
\end{multline}
for $0 < \rho < R <1$. Since for some $w, \;p$ the exponents $2w - p +2, 2-p$ may be nonpositive, 
we estimate respective powers of $ |\nabla u| $ using $ |\nabla u|^s $ as follows
\begin{multline}
\int_{Q_R} | \nabla u|^s  [(1+s) ( | \nabla u|^{2w - p +2} + | \nabla u|^{2-p} )] dxdt
 \le (1+s) \int_{Q_R} [1 + | \nabla u|^{\max(s+2 -p, s+2w +2 -p)}] dxdt
 \label{bg_p11}
 \end{multline}
 For the last inequality to hold, we must assume 
 \begin{equation}
 s > \max (p-2w-2,  \: p-2  )
  \label{bg_p12}
 \end{equation}
 Because summand $ | \nabla u|^{2-p}$ occurs only if  $c_2 \neq 0$ in the growth condition $(A1)$: $ |f^i| \le  c_1 | \nabla u |^w + c_2$,  the above assumption   \eqref{bg_p12} can be written as
 \begin{equation*}
 \begin{cases}
s  > p-2 \; for\;  c_2 \neq 0 \\
s  > p-2w-2 \; for\;  c_2 = 0
\end{cases}
  \label{bg_p13}
\end{equation*}
In the forthcoming iteration scheme we construct a growing sequence of $s_i$, therefore it is sufficient to assume 
 \begin{equation*}
 \begin{cases}
s_0  > p-2 \; for\;  c_2 \neq 0 \\
s_0  > p-2w-2 \; for\;  c_2 = 0
\end{cases}
  \label{bg_p14}
\end{equation*}
which coincides with our initial integrability assumption $(A2)$.
\newline
By computation, the following inequality holds
  \begin{multline}
  \int_{Q_R} | \nabla( | \nabla u|^{\frac{p+s}{2}} \eta )|^2 dxdt \le 
   \int_{Q_R} | \nabla u|^{p+s} |\nabla \eta|^2 dxdt + 
  (p+s)^2 \int_{Q_R}
  |\nabla u|^{p+s-2} |\nabla^2 u|^2 \eta^2 dxdt
   \label{bg_p8}
  \end{multline}
 \noindent
Adding to both sides $\frac{ \lambda - \varepsilon}{(p+s)^2} \int_{Q_R} | \nabla u|^{p+s} |\nabla \eta|^2 dxdt$ and considering properties of $\eta$, as $(s+p)^2 \le c (1+s^2)$, we arrive from \eqref{bg_p10},  by virtue of \eqref{bg_p8},  at
\begin{equation}
\sup_t \int_{B_R} |\nabla u|^{s+2} \eta^2 dx + \int_{Q_R} | \nabla( | \nabla u|^{\frac{p+s}{2}} \eta )|^2 dxdt \le
 \frac{c (1+ s^3)}{(R -\rho)^{M} }\int_{Q_R} 1 + | \nabla u|^{s+M} dxdt
  \label{bg_p15}
 \end{equation}
taking into account  \eqref{bg_p11} if neccessary. Recall that by definition $M = \max (2, p, 2q-p, w+1, 2w-p+2)  $. \\
By H\"older and critical-Sobolev inequalities (respectively), one gets
\begin{multline}
  \int_{Q_\rho} | \nabla u|^{p+s +(s+2) \frac{2}{n}} dxdt \le
  \int_{t_0 - R^2}^{t_0} \left[ \int_{B_R} |\nabla u|^{s+2} \eta^2 dx \right]^{\frac{2}{n}} \left[  \int_{B_R}  | \nabla u|^{(s+p) \frac{n}{n-2}} \eta^{ \frac{2 n}{n-2}}  dx \right]^{\frac{n-2}{n}} dt \cr
\le     \left[ \sup_t \int_{B_R} | \nabla u|^{s+2} \eta^2 dx \right]^{\frac{2}{n}}  \int_{t_0 - R^2}^{t_0} \left[  \int_{B_R}  (| \nabla u|^{\frac{p+s}{2}} \eta)^{ \frac{2n}{n-2}} dx \right]^{\frac{n-2}{n}} dt \le \\
    \left[ \sup_t \int_{B_R} |\nabla u|^{s+2} \eta^2 dx \right]^{\frac{2}{n}} \int_{Q_R}  | \nabla (| \nabla u|^{ \frac{p+s}{2} } \eta)|^2 dxdt \stackrel{\eqref{bg_p15}}{\le}
\left[ \frac{c (1+s^3)}{ (R - \rho)^{M}}  \int_{Q_R} 1+ |\nabla u|^{s+M} dxdt
 \right]^{1+\frac{2}{n}}
\label{bg_energ}
  \end{multline}
Inequality \eqref{bg_energ} is our desired energy estimate, which we now iterate.
Define recursively numbers $s_i$: $s_{i+1} + M = p+s_i +(s_i+2) \frac{2}{n}$, then
\begin{equation*}
s_{i} = \left(1+\frac{2}{n} \right)^i \left[s_0 + n+2 - \frac{ n (p - M)}{2} \right] - \left[ 2 - \frac{ n (p- M)}{2} \right]
\label{bs_n>p_rek}
\end{equation*}
Utilizing the initial-integrability assumption, i.e. $s_0 + 2+ n - \frac{n (p- M)}{2} > 0$, we have
\begin{equation}
s_{i} \xrightarrow{i\to\infty} \infty; \quad \frac{s_{i}}{(1+\frac{2}{n})^i} \xrightarrow{i\to\infty}
s_0 +2 +n - \frac{n (p-M)}{2}
\label{bg_n>p_rek1}
\end{equation}
Let
\begin{equation*}
\psi_i = \dashint_{S_{R_i}} |\nabla u|^{s_i +M} dxdt
\end{equation*}
then \eqref{bg_energ} with $\eta_{R_{i+1}, R_i} $ can be written as
\begin{multline}
|S_{R_{i+1}}| \psi_{i+1} \le \left[C (1+ s_i^p) \left(\frac{2^{i+2}}{R_0}\right)^M  |S_{R_{i}}| (1+ \psi_i ) \right]^{1+ \frac{p}{n}} \implies \cr 
R^{n+2}_{i+1} \psi_{i+1} \le \left[C (1+ s_i^p) \left(\frac{2^{i+2}}{R_0} \right)^M  R^{n+2}_{i+1} (1+ \psi_i ) \right]^\beta \implies
\psi_{i+1} \le C (1+ s_i^p)^\beta 2^i (1+ \psi_i ) ]^\beta
\label{bg_n>p_rek2}
\end{multline}
with $ \beta: = 1+ \frac{2}{n}$; the last inequality given by  $R_i := \frac{R_0}{2} (1+ 2^{-i})$. As we know from \eqref{bg_n>p_rek1} that asymptotically $s_i$ behaves like $\beta^i$, finally \eqref{bg_n>p_rek2} folds to
\begin{equation*}
\psi_{i+1} \le C^i \psi_i^\beta + C^i
\label{bg_n>p_rek3}
\end{equation*}
which, by a standard computation (see \citet{Bae-Choe04} for details), gives
\begin{equation}
\psi_{i+1} \le C^{\beta^{i+1}} \psi_0^{\beta^{i+1}}  + (i+1) C^{\beta^{i+1}}
\label{bg_n>p_rek4}
\end{equation}
From the above considerations one gets, using the definition of $\psi$
\begin{multline}
R_0^{ - \frac{n+2}{s_{i+1} +M }}    | \nabla u|_{L^{s_{i+1} + M} \left(Q_{\frac{R_0}{2}}\right)} \le \left( \dashint_{S_{R_{i+1}}} |u|^{s_{i+1} +M} \right)^{\frac{1}{s_{i+1} +M }}=  \psi_{i+1}^{\frac{1}{s_{i+1} +M }} \stackrel{\eqref{bg_n>p_rek4}}{\le} \cr  (C^{\beta^{i+1}} \psi_0^{\beta^{i+1}})^{\frac{1}{s_{i+1} +M }}  + ((i+1) C^{\beta^{i+1}})^{\frac{1}{s_{i+1} +M }}\xrightarrow{i\to\infty} C \psi_0^{\frac{1}{s_{0} +2 +n - \frac{Mn}{p} }} + C
\label{bg_n>p_rek5}
 \end{multline}
  in view of \eqref{bg_n>p_rek1}. \\
  As $s_i + M \xrightarrow{i\to\infty} \infty$, \eqref{bg_n>p_rek5} in tandem with the initial integrability assumption gives the following uniform bound
\begin{equation*}
 |\nabla u|_{L^{\infty} (Q_{\frac{R_0}{2}})} \le C \left( \dashint_{Q_{R_{0}}} |\nabla u|^{s_0 +M} dxdt \right)^{\frac{1}{s_{0} +2 +n - \frac{Mn}{p} }} + C
 \end{equation*}

\end{proof}
\noindent
In the next theorem we resign from the term possessing $q>p$ in the ellipticity assumption. This allows us, in turn, to obtain bigger growths of the right-hand-side, as now it is possible to derive estimates for negative $s > - \frac{\lambda}{\Lambda}$.
\begin{Tw}
Gradient of solution to \eqref{bg1} is locally bounded, under the following assumptions:
\begin{enumerate}
\item[(A0)] ellipticity-type: $A^i_{\alpha}$ is given by potential $F \in C^2 (\er)$ as follows
\begin{equation*}
A^i_\alpha (Q) = (F(|Q|)_{Q^i_\alpha}
\end{equation*}
and $F$ enjoys ellipticity
\begin{equation*}
(F(|Q|)_{Q^i_\alpha Q^j_\beta} \zeta^i_\alpha \zeta^j_\beta  \in [ \lambda  | Q|^{p-2} ,  \Lambda | Q|^{p-2}] |\zeta|^2
\end{equation*}
\item[(A1)] growth-type
\begin{equation*}
|f^i| \le c | \nabla u |^w, \quad p \ge w \ge 0, \quad c \in L^\infty (\Omega_T)
\end{equation*}
\item[(A2)] initial integrability
\begin{equation*}
\nabla u \in L_{loc}^{s_0 + M} (\Omega_T), \quad M := \max (2, p, w+1, 2w-p+2)  
\end{equation*}
with $s_0 $ satisfying
\begin{equation*}
\begin{cases}
s_0 > \max (- \frac{\lambda}{\Lambda}, p-2w-2) \\
s_0 + 2 +\frac{np}{2} - \frac{M n}{2}  > 0 
\end{cases}
\end{equation*}
\end{enumerate}
moreover, for any $Q_{R_{0}} \subset \Omega$ with $R_{0} < 1$ following inequality holds
\begin{equation*}
 |\nabla u|_{L^{\infty} (Q_{\frac{R_0}{2}})} \le C \left( \dashint_{Q_{R_{0}}} |\nabla u|^{s_0 +M} dxdt \right)^{\frac{1}{s_0 +2 + \frac{np}{2} - \frac{Mn}{2}}} + C
 \end{equation*}
\label{thm2_p=q}
\end{Tw}
\noindent
\begin{proof} For $s \ge 0$ Theorem \ref{thm2_p=q} is a special case of  Theorem~\ref{thm2}, therefore it suffices to show it in the case of negative $s$. The only difference in the energy estimates is the lack of positivity of the left-hand-side term of  \eqref{bg_p2}, where sign of $s$ plays a role
\begin{equation*}
\int_{B_R} A^i_{\alpha \: '(\nabla u)^j_\beta} (\nabla u)   u^j_{x_\beta  x_\gamma}  s |\nabla u|^{s-2} u^i_{x_\gamma} u^k_{x_\gamma x_\delta}   u^k_{x_\delta} \eta^2 dx
\end{equation*}
it can be however estimated as follows
 \begin{multline*}
 \int_{B_R} A^i_{\alpha \: '(\nabla u)^j_\beta} (\nabla u)   u^j_{x_\beta  x_\gamma} s |\nabla u|^{s-2} u^i_{x_\gamma} u^k_{x_\gamma x_\delta}   u^k_{x_\delta} \eta^2 dx =\\  s \int_{B_R}  \eta^2 F''(|\nabla u |) {|\nabla u |^{s-4}}  u^i_{x_\gamma}  u^j_{x_\beta} u^j_{x_\beta x_\gamma}  u^i_{x_\alpha}  u^k_{x_\delta} u^k_{x_\delta x_\alpha} dx \ge
 s \Lambda \int_{B_R}  \eta^2 | \nabla u |^{p+s-2} | \nabla^2 u|^2 dx
\end{multline*}
which allows for a following counterparty of  \eqref{bg_p6}
\begin{multline*}
 \int_{B_R} A^i_{\alpha \: '(\nabla u)^j_\beta} (\nabla u)   u^j_{x_\beta  x_\gamma} ( |\nabla u|^s u^i_{x_\gamma x_\alpha} \eta^2 + s |\nabla u|^{s-2} u^i_{x_\gamma} u^k_{x_\gamma x_\delta}   u^k_{x_\delta} \eta^2 ) dx \\
 \ge  (\lambda + s \Lambda) \int_{B_R} |\nabla u|^{p-2}   |\nabla^2 u|^2  |\nabla u|^{s}  \eta^2 dx
\end{multline*}
From this inequality on, one proceeds identically as in proof of Theorem \ref{thm2}.
\end{proof}
\noindent
As a consequence of  Theorem \ref{thm2_p=q} we obtain Theorem \ref{thm_intr}.
\begin{proof}[Proof of Theorem \ref{thm_intr}.]
The point $1.$ stems from theory in \citet{diBbook93} and it implies that in points $2.$ and $3.$ one can assume $w >p-1$. The rest of the proof follows from Theorem \ref{thm2_p=q}. Indeed, for p-Laplace system $\lambda = \Lambda$, so the ellipticity condition is fulfilled and the rest of the assumptions of Theorem \ref{thm2_p=q} can be rewritten with $s_0 = \tilde{p} - M $ as follows:
\begin{equation*}
M := \max (2, 2w-p+2)  
\end{equation*}
\begin{equation*}
\nabla u \in L_{loc}^{\tilde{p}} (\Omega_T)
\end{equation*}
\begin{equation}
 \tilde{p} > \max \left( M -1, \frac{5M}{2} - \frac{3p}{2} - 2 \right)
 \label{tilde}
\end{equation}
In the case $w \ge \frac{p}{2}$: $M = 2w-p+2 $ and  \eqref{tilde} gives $w < \frac{ {\tilde p} + p - 1}{2}$ for $w \le p- \frac{2}{3}$ and $w < \frac{ {\tilde p} + 4p - 3}{5}$ for $w \ge p- \frac{2}{3}$; as in view of  $w >p-1$ and $\tilde{p} \ge p$ the first condition is void, we obtain assumption $w < \frac{ {\tilde p} + 4p - 3}{5}$.\\
In the case $w \le \frac{p}{2}$: $M=2$ and  \eqref{tilde}: for $p \ge \frac{4}{3} $ takes the form  $\tilde{p} > 1$, which always holds, and for $p < \frac{4}{3} $ it reads  ${\tilde p} > 3 - \frac{3p}{2}$. These two conditions are: $ p > 2 - \frac{2}{3} \tilde {p}$.
\end{proof}
\noindent
Please recall, that  by point $2.$ merely from existence, i.e for ${\tilde p} = p$, one has $w < p - \frac{3}{5}$.

\section{Conclusion.}
\subsection{Note on the rigorous estimates.}
\noindent
The above computations are formal. To perform them rigorously, transform the considered problem using Steklov averages with respect to time and use finite differences instead of differentiating it with respect to space. This procedure has been presented in \citet{diBbook93, diB-Frie84} for homogenous systems. In our case we need to deal additionally with a quasilinear right-hand-side, for which the testing function  $u^i_{x_\gamma} | \nabla u|^s \eta^2$ may not be admissible. In order to begin iteration, we need to have: $s_0 +w+1 \le p$ and to perform it at the i-th step: $s_i +w +1 \le s_i + M$. However the latter inequality holds from the definition of $M$, the former may prove sometimes troublesome. In such cases one can resort to testing with  $F_n (u^i_{x_\gamma} | \nabla u|^s \eta^2)$, where $F_n (x) $ is a Lipschitz truncation at the level $n$. As the estimates are valid for every $n$ we can proceed as before. Observe however, that we do not encounter these difficulties during computations for Theorem \ref{thm_intr}. For additional rigorous treatment (especially for $s $ nonpositive, consult \citet{Bae-Choe04} and reference therein\footnote{Observe, however, that in \citet{Bae-Choe04} there are allowed $s_0 > -2$, which seems to be incorrect as far as $s \in (-2, -1]$ are concerned.} as well as \citet{choe92, choe91}). 
\subsection{Further research.}
\noindent
The most obvious possible generalizations to the result is to allow for faster growths of the right-hand-side at the cost of additional assumptions, especially as some of them, like boundedness of the solution, appear naturally in the existence theory. It would be interesting to obtain a general result for critically growing right-hand-side (i.e. like $1+ |\nabla u|^p$), with some smallness assumption, which would generalize the classical results for the heat system, mentioned in the introduction.




\bibliographystyle{plainnat}
\bibliography{ogr}



\end{document}